\documentclass[twoside,10pt,reqno]{amsart}
\usepackage{amsmath, amsthm,amssymb,amstext,amsfonts,amscd}
\usepackage{bigints}

\usepackage{graphicx}
\usepackage{caption}
\usepackage{refstyle}
\usepackage{multirow}
\usepackage{hyperref}
\usepackage{lmodern}
\usepackage{latexsym}
\usepackage{amssymb}
\usepackage[utf8]{inputenc}
\usepackage[english]{babel}
\usepackage{wasysym} 
\usepackage{cite}
\numberwithin{equation}{section}
\begin{document}

\begin{center}
\textbf{ANALYTICAL EXPRESSIONS FOR THE EXACT CURVED SURFACE AREA OF A HEMIELLIPSOID VIA MELLIN-BARNES TYPE CONTOUR INTEGRATION}
\end{center}
\begin{center}
\noindent {M.A. Pathan$^{1,2}$, M. I. Qureshi$^3$, Javid Majid$^{3,*}$    }
\end{center}
\begin{center}
$^1$ Centre for Mathematical and Statistical Sciences (CMSS), Peechi,\\
Thrissur, Kerala-680653, India\\
$^2$ Department of Mathematics, Aligarh Muslim University,\\
Aligarh, U.P., India\\
$^3$ Department of Applied Sciences and Humanities \\
Faculty of Engineering and Technology\\
Jamia Millia Islamia (A Central University), New Delhi-110025, India.\\
 Emails: mapathan@gmail.com,~miqureshi\_delhi@yahoo.co.in\\
  $^*$Corresponding author: javidmajid375@gmail.com 
\end{center}
\textbf{Abstract:} 
In this article, we aim at obtaining the analytical expressions ({\bf not previously found and not recorded in the literature}) for the exact curved surface area of a hemiellpsoid in terms of Appell's double hypergeometric function of first kind. The derivation is based on Mellin-Barnes type contour integral representations of generalized hypergeometric function$~_pF_q(z)$, Meijer's $G$-function and analytic continuation formula for Gauss function. Moreover, we obtain some special cases related to ellipsoid, Prolate spheroid and Oblate spheroid. The closed forms for the exact curved surface area of a hemiellpsoid are also verified numerically by using {\it Mathematica Program}.\\



\section{\bf{Introduction and preliminaries}}
\noindent
For the definition of Pochhammer symbols, power series form of generalized hypergeometric function ${}_pF_q(z)$ and several related results, we refer the beautiful monographs (see, e.g., \cite{And1999, Erd1, Lebedev1972,Luke169 , Rain, Slater1966, SriMano})\\

\noindent
$\bullet$ Some results recorded in the table of Prudnikov {\em et al.}[\cite{Prudbrimari90}, p.474, Entry(98) and p.479, Entry(210)]:
\begin{equation}\label{Z1.9}
~_2{F}_1\left[\begin{array}{lll}
\frac{1}{2},~2;\\
& z\\
~~~\frac{3}{2};\\\end{array}\right]=\frac{1}{2}\left[\frac{1}{(1-z)}+\frac{\tanh^{-1}(\sqrt{z})}{\sqrt{(z)}}\right];~~|z|<1,~~~~~~~~~~~~ 
\end{equation}
\begin{equation}\label{Z1.8}
~_2{F}_1\left[\begin{array}{lll}
\frac{3}{2},2;\\
& z\\
~~~\frac{5}{2};\\\end{array}\right]=\frac{3}{2z}\left[\frac{1}{(1-z)}-\frac{\tanh^{-1}(\sqrt{z})}{\sqrt{z}}\right];~~|z|<1.~~~~~~~~~~~~
\end{equation}

\noindent
$\bullet$  Analytic continuation formula \cite[p.63, Eq.(2.1.4(17)), \cite{Lebedev1972}, p.249, Eq.(9.5.9),\cite{Prudbrimari90}, p.454, Entry(7.3.1(6)), \cite{Slater1966}, p.36, Eq.(1.8.1.11)]{Erd1}:\\

\noindent
When$|z|>1$, then
\begin{equation*}
~\mathbf~{_{2}F_{1}}\left[\begin{array}{ccc}
~a,~b;   \\
& z\\
~~~~c; \\
\end{array}\right]
=\frac{\Gamma(c)~\Gamma(b-a)}{\Gamma(b)~\Gamma(c-a)}(-z)^{-a}~{_{2}F_{1}}\left[\begin{array}{ccc}
a,1+a-c;   \\
& \frac{1}{z}\\
~~~1+a-b; \\
\end{array}\right]+~~~~~~~~~~~~~~~~~~~~~~
\end{equation*}
\begin{equation}\label{G7.12}
+\frac{\Gamma(c)~\Gamma(a-b)}{\Gamma(a)~\Gamma(c-b)}(-z)^{-b}~{_{2}F_{1}}\left[\begin{array}{ccc}
b,1+b-c;   \\
& \frac{1}{z}\\
~~~1+b-a; \\
\end{array}\right],~~~~~~~~~~~~~~~~~~~~~~
\end{equation}
where $|\arg(-z)|<\pi,~|\arg(1-z)|<\pi$ and $(a-b)\neq 0,\pm1,\pm2,\pm3,....$\\

\noindent

\noindent
$\bullet$ Mellin-Barnes type contour integral representation of binomial function:
\begin{equation}\label{D12.100}
(1-z)^{-a}=~{}_1F_0\left[\begin{array}{ll}
~a;\\
&z\\
-;\end{array}\right]=\frac{1}{(2\pi i)~\Gamma(a)}\int_{-i\infty}^{+i\infty}\Gamma(a+s)\Gamma(-s)(-z)^s~ds:~~z\neq 0,
\end{equation}
where $|\arg(-z)|<\pi,|z|<1,a\in\mathbb{C}\backslash\mathbb{Z}^-_0$ and $i=\sqrt{(-1)}$.\\

\noindent
$\bullet$ Appell's function of first kind \cite[p.53, Eq.(4)]{SriMano} is defined as:
\begin{equation}\label{C1.4}
F_1\left[\begin{array}{ll}
a;~b,c;~d;~x,y\end{array}\right]=F^{1:1;1}_{1:0;0}\left[\begin{array}{ll}
~a:~b;~c;\\
& x,y\\
~d:-;-;\\\end{array}\right]=\sum_{m,n=0}^{\infty}\frac{(a)_{m+n}(b)_m(c)_n~x^m~y^n}{(d)_{m+n} ~ m! ~n!}
\end{equation}
\begin{equation}\label{Y1.4}
=\sum_{m=0}^{\infty}\frac{(a)_{m}(b)_m~x^m}{(d)_{m} ~ m! }~{_{2}F_{1}}\left[\begin{array}{ccc}
                                        a+m,~c;   \\
                                        & y\\
                                      ~~  d+m; \\
                                        \end{array}\right]=\sum_{n=0}^{\infty}\frac{(a)_{n}(c)_n~y^n}{(d)_{n} ~ n! }~{_{2}F_{1}}\left[\begin{array}{ccc}
                                                                                a+n,~b;   \\
                                                                                & x\\
                                                                              ~~  d+n; \\
                                                                                \end{array}\right]
\end{equation}
$\bullet$ Convergence conditions of Appell's double series $F_1$:
\begin{enumerate}
\vsize.1cm
\item[{(i)}] Appell's series $F_1$ is convergent when $|x|<1,|y|<1;~a,b,c,d\in\mathbb{C}\backslash\mathbb{Z}^-_0$.

\item[{(ii)}] Appell's series $F_1$ is absolutely convergent when $|x|=1,|y|=1;~a,b,c,d\in\mathbb{C}\backslash\mathbb{Z}^-_0;~\mathfrak{R}(a+b-d)<0,\mathfrak{R}(a+c-d)<0$ and $\mathfrak{R}(a+b+c-d)<0$.

\item[{(iii)}] Appell's series $F_1$ is conditionally convergent when $|x|=1,|y|=1;~x\neq 1,y\neq 1;~a,b,c,d\in\mathbb{C}\backslash\mathbb{Z}^-_0;~\mathfrak{R}(a+b-d)<1,\mathfrak{R}(a+c-d)<1$ and $\mathfrak{R}(a+b+c-d)<2$.

\item[{(iv)}] Appell's series $F_1$ is a polynomial If $a$ is a negative integer;  $b,c,d\in\mathbb{C}\backslash\mathbb{Z}^-_0$.

\item[{(v)}] Appell's series $F_1$ is a polynomial If $b$ and $c$ are negative integers;  $a,d\in\mathbb{C}\backslash\mathbb{Z}^-_0$.
\end{enumerate}

\noindent
$\bullet$ Mellin-Barnes type contour integral representation of Meijer's G-function (\cite[p.45, Eq.(1)]{SriMano}, see also\cite{Erd1,Luke169}):\\
When $p\leq q$ and $1\leq m \leq q,~0\leq n\leq p,$ then
\begin{equation*}
G^{m,n}_{p,q}\left(z\left|   \begin{array}{ll}
\alpha_1,\alpha_2,\alpha_3,...,\alpha_n;\alpha_{n+1},...,\alpha_p\\
\beta_1,\beta_2,\beta_3,...,\beta_m;\beta_{m+1},...,\beta_q  \end{array}\right. \right)=\frac{1}{2\pi i}\int_{-i\infty}^{+i\infty}\frac{\prod_{j=1}^{m}\Gamma(\beta_j-s)\prod_{j=1}^{n}\Gamma(1-\alpha_j+s)}{\prod_{j=m+1}^{q}\Gamma(1-\beta_j+s)\prod_{j=n+1}^{p}\Gamma(\alpha_j-s)}(z)^s~ds
\end{equation*}
\begin{equation}\label{E12.12}
=\frac{1}{2\pi i}\int_{-i\infty}^{+i\infty}\frac{\Gamma(\beta_1-s)...\Gamma(\beta_m-s)\Gamma(1-\alpha_1+s)...\Gamma(1-\alpha_n+s)}{\Gamma(1-\beta_{m+1}+s)...\Gamma(1-\beta_q+s)\Gamma(\alpha_{n+1}-s)...\Gamma(\alpha_p-s)}(z)^s~ds,
\end{equation}
where $z\neq 0,(\alpha_i-\beta_j)\neq$ positive integers, $i=1,2,3,...,n;~j=1,2,3,...,m$. For details of contours, see\cite[p.207, \cite{Luke169}, p.144]{Erd1}.\\

\noindent
$\bullet$ Convergence conditions of Meijer's G-function:\\
When $\Lambda=m+n-\left(\frac{p+q}{2}\right),~\nu=\sum_{j=1}^{q}\beta_j-\sum_{j=1}^{p}\alpha_j, $ then
\begin{enumerate}
\vsize.1cm
\item[{(i)}] The integral (\ref{E12.12}) is convergent when $|arg(z)|<\Lambda\pi$ and $\Lambda> 0$.\\

\item[{(ii)}] If $|arg(z)|=\Lambda\pi$ and $~\Lambda\geq 0$, then the integral (\ref{E12.12}) is absolutely convergent when $p=q$ and $\mathfrak{R}(\nu)<-1$.

\item[{(iii)}]If $|arg(z)|=\Lambda\pi$ and $~\Lambda\geq 0$, then the integral (\ref{E12.12}) is also absolutely convergent, when $p\neq q,~(q-p)\sigma>\mathfrak{R}(\nu)+1-\left(\frac{q-p}{2}\right)$ and $s=\sigma+ik$, where $\sigma$ and $k$ are real. $\sigma$ is chosen so that for $k  \rightarrow\pm\infty$.\\

\noindent
For other two types of contours, following will be convergence conditions of the integral (\ref{E12.12}):\\

\item[{(iv)}]The integral (\ref{E12.12}) is  convergent if $q\geq 1$ and either $p<q,0<|z|<\infty$ or $p=q,~0<|z|<1$.

\item[{(v)}]The integral (\ref{E12.12}) is  convergent if $p \geq 1$ and either $p>q,0<|z|<\infty$ or $p=q,~|z|>1$.\\
\end{enumerate}

\noindent
$\bullet$ Relations between Meijer's $G$- function and ${}_2F_1(z)$ \cite[p.61, \cite{Wikipedia}, p.77, Eq.(1)]{Mathaisaxena}: 
\begin{equation*}
G^{2~2}_{2~2}\left(z\left|   \begin{array}{ll}
1-a,1-b;-\\
0,c-a-b;-  \end{array}\right. \right)=\frac{\Gamma(a)\Gamma(b)\Gamma(c-a)\Gamma(c-b)}{\Gamma(c)}{}_2F_1\left[\begin{array}{ll}
a,~b;\\
& 1-z\\
~~~ ~c;\end{array}\right],
\end{equation*}
where $|1-z|<1$ and $c-a,c-b\neq 0,-1,-2,...$
\begin{equation*}
G^{2~2}_{2~2}\left(z\left|   \begin{array}{ll}
a_1,a_2;-\\
b_1,b_2;-  \end{array}\right. \right)=\frac{\Gamma(1-a_1+b_1)\Gamma(1-a_1+b_2)\Gamma(1-a_2+b_1)\Gamma(1-a_2+b_2)z^{b_1}}{\Gamma(2-a_1-a_2+b_1+b_2)}\times
\end{equation*}
\begin{equation}\label{F12.200}
\times{}_2F_1\left[\begin{array}{ll}
1-a_1+b_1,1-a_2+b_1;\\
& 1-z\\
~~~2-a_1-a_2+b_1+b_2;\end{array}\right];~~~|1-z|<1.
\end{equation}

\noindent 
$\bullet$ When the two dimensional curve i.e, generating curve lying in $x$-$y$ plane (suppose $y=f(x)$) is revolved about $x$-axis, then equation of generated three dimensional surface will be $y^2+z^2=[f(x)]^2$.\\
 \noindent
$\bullet$ When the two dimensional curve i.e, generating curve lying in $x$-$y$ plane (suppose $x=F(y)$) is revolved about $y$-axis, then equation of generated three dimensional surface will be $x^2+z^2=[F(y)]^2$.\\
\noindent
Similarly we can write the equation of generated three dimensional surface, when the curve lies in $y$-$z$ plane and $z$-$x$ plane.\\
  
\noindent 
$\bullet$ The equation of an ellipse is
\begin{equation}\label{l12.70}
\frac{x^2}{a^2}+\frac{y^2}{b^2}=1;~~a> b>0.~~~~~~~~~~~~~~~
\end{equation}
$\bullet$ When the above ellipse  represented by (\ref{l12.70}) is revolved about major axis (i.e, $x$-axis), then equation of generated 3-D surface called Prolate spheroid, will be

\begin{equation}\label{l12.71}
\frac{x^2}{a^2}+\frac{y^2+z^2}{b^2}=1;~~a> b>0.~~~~~~~~~~~~~~~
\end{equation}
\begin{figure}[ht!]
\noindent\begin{minipage}{0.45\textwidth}
\includegraphics[width=\linewidth]{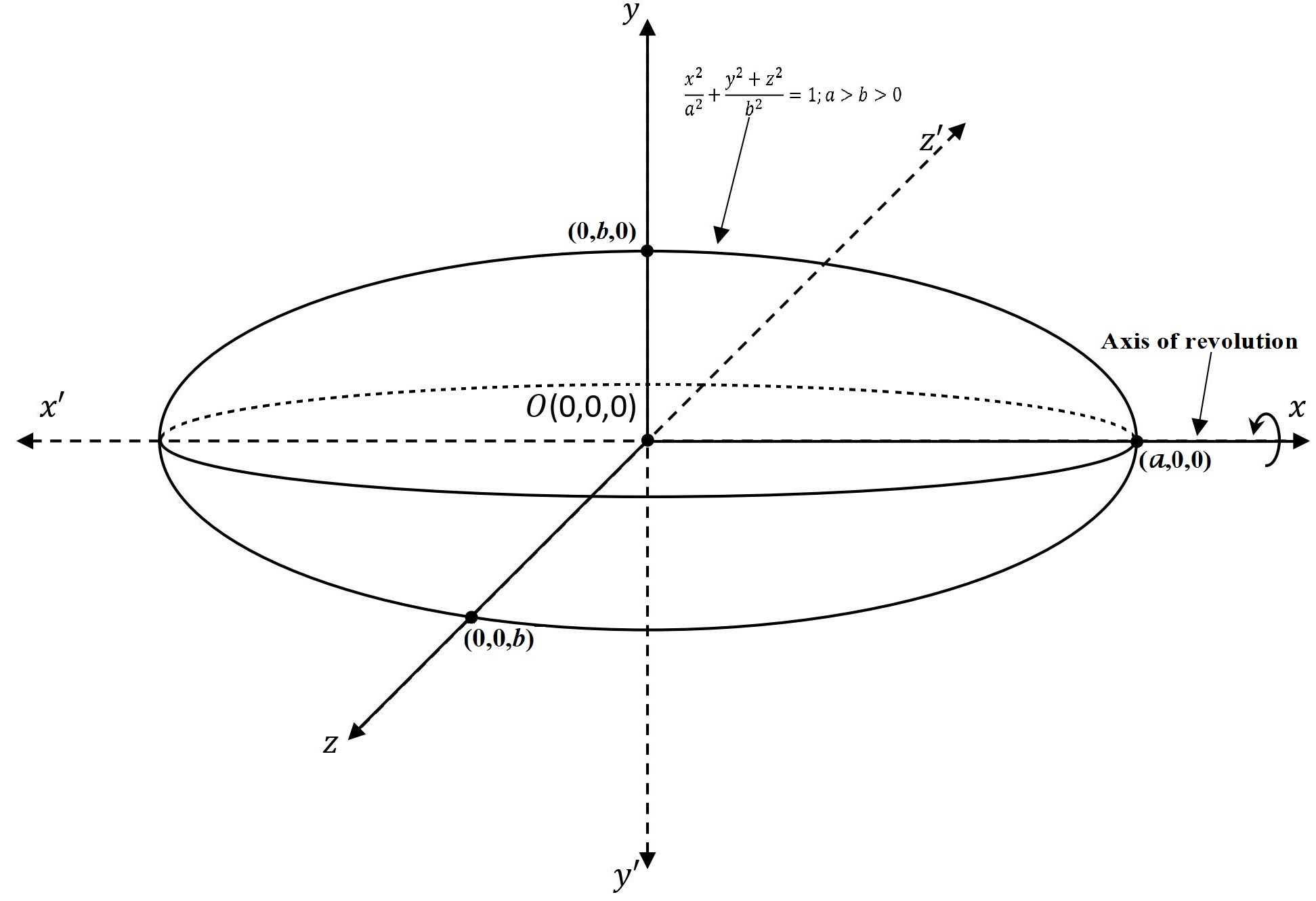}
\end{minipage}%
\caption[short]{ Prolate Spheroid.}
\label{fig:Projection }
\end{figure}

\noindent
$\bullet$ When the  ellipse represented by (\ref{l12.70}) is revolved about minor axis (i.e, $y$-axis), then equation of generated 3-D surface called Oblate spheroid, will be
\begin{figure}[ht!]
\noindent\begin{minipage}{0.50\textwidth}
\includegraphics[width=\linewidth]{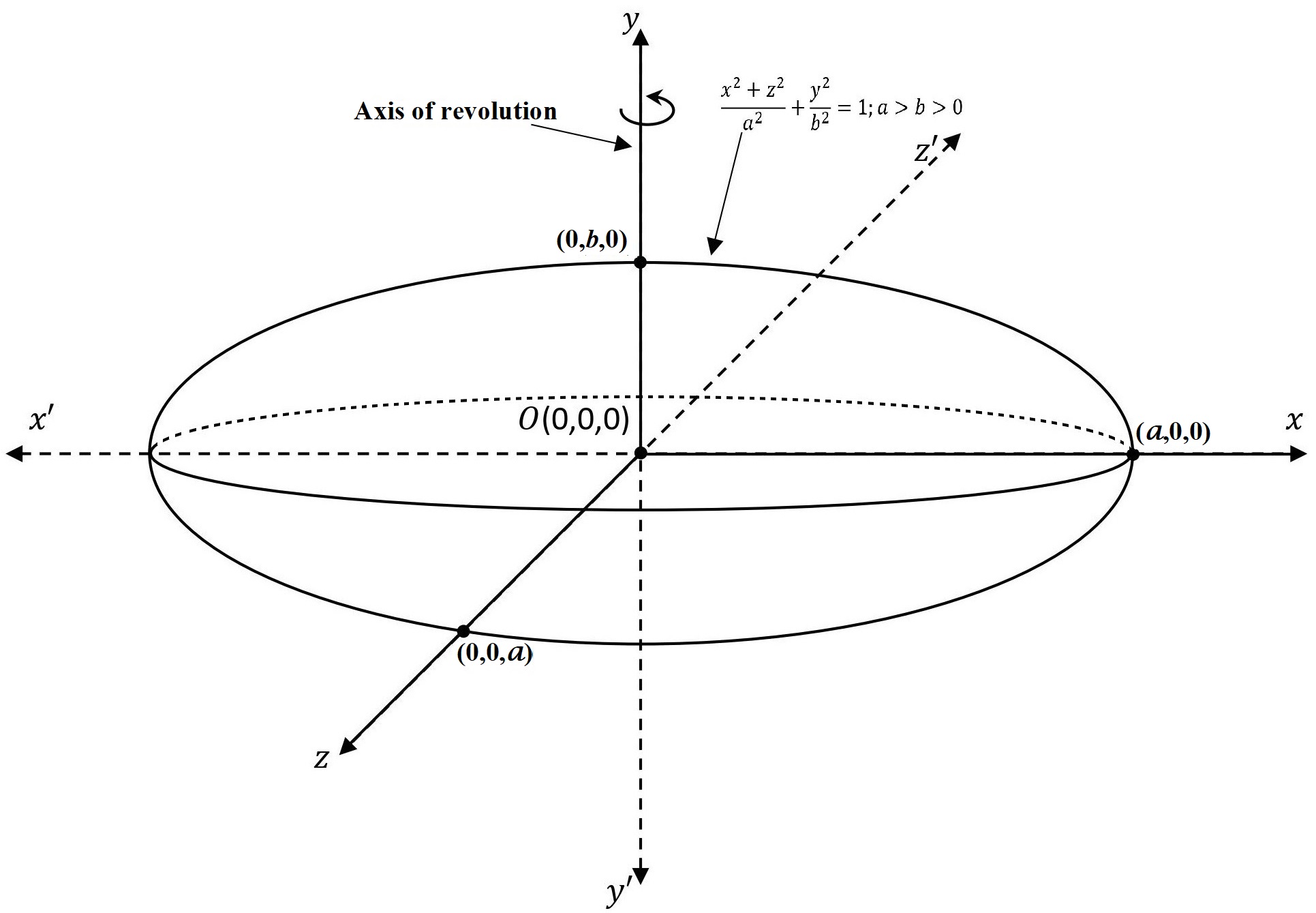}
\end{minipage}%
\caption[short]{Oblate Spheroid.}
\label{fig:Projection }
\end{figure}

\begin{equation}\label{l12.72}
\frac{x^2+z^2}{a^2}+\frac{y^2}{b^2}=1;~~a> b>0.~~~~~~~~~~~~~~~
\end{equation}

\noindent
$\bullet$ The sphere $x^2+y^2+z^2=c^2$, prolate spheroid and oblate spheroid are the particular cases of the ellipsoid 
\begin{equation}
\frac{x^2}{a^2}+\frac{y^2}{b^2}+\frac{z^2}{c^2}=1;~~a> b>0.~~~~~~~~~~~~~~~
\end{equation}

\noindent
$\bullet$ Suppose $\phi(x,y)=0$ is the projection of the curved surface of three dimensional figure $z=f(x,y)$ over the $x$-$y$ plane, then curved surface area is given by
\begin{equation}\label{a12.9}
\hat{S}=\underbrace{\int\int}_{\stackrel{\text{over the area }} {\phi(x,y)=0}}\sqrt{\left\lbrace 1+\left( \frac{\partial z}{\partial x}\right) ^2+\left( \frac{\partial z}{\partial y}\right) ^2\right\rbrace }~dx~dy.
\end{equation}
$\bullet$ Suppose $\psi(y,z)=0$ is the projection of the curved surface of three dimensional figure $x=g(y,z)$ over the $y$-$z$ plane, then curved surface area is given by
\begin{equation}\label{t12.9}
\hat{S}=\underbrace{\int\int}_{\stackrel{\text{over the area }} {\psi(y,z)=0}}\sqrt{\left\lbrace 1+\left( \frac{\partial x}{\partial y}\right) ^2+\left( \frac{\partial x}{\partial z}\right) ^2\right\rbrace }~dy~dz.
\end{equation}
\begin{figure}[ht!]
\noindent\begin{minipage}{0.45\textwidth}
\includegraphics[width=\linewidth]{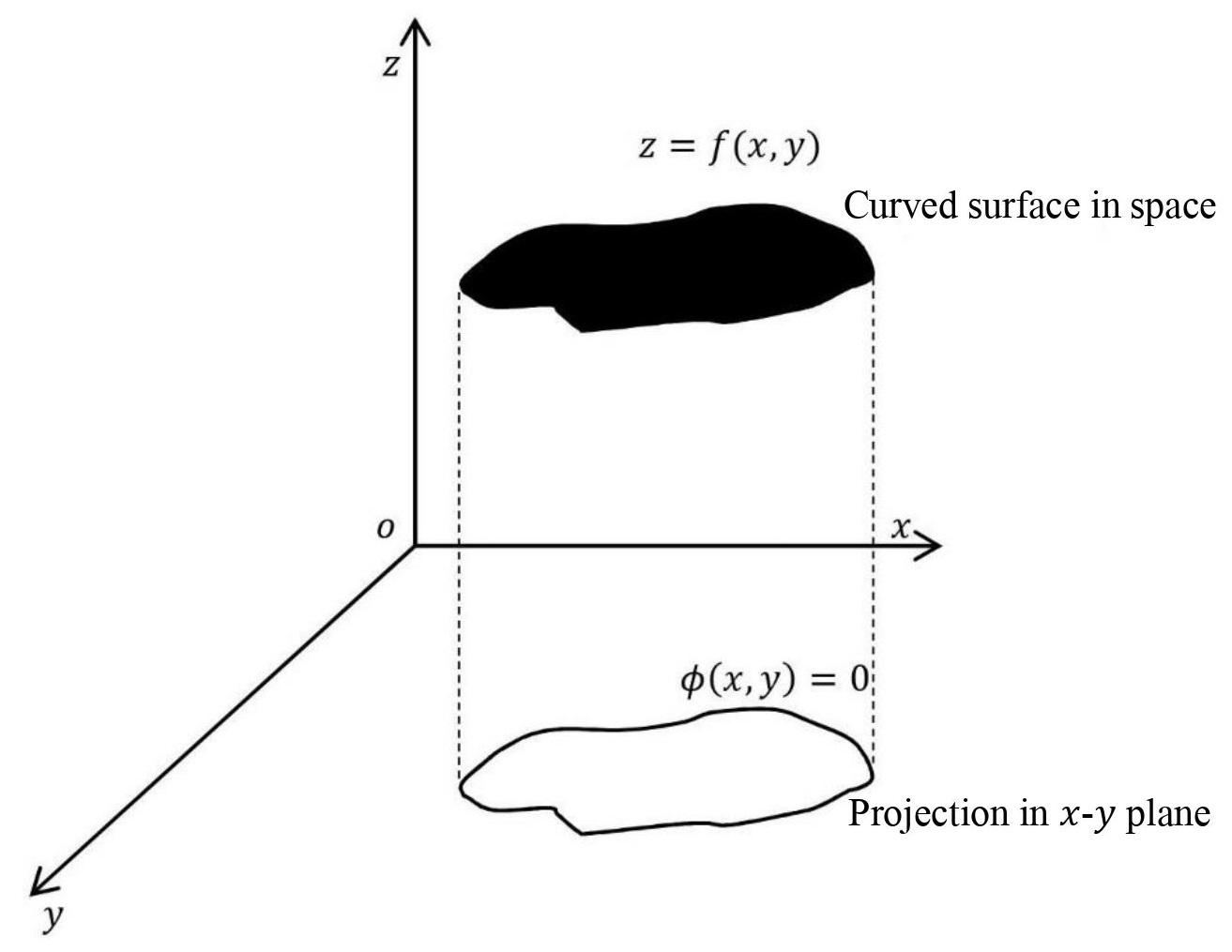}
\end{minipage}%
\caption[short]{Projection of curved surface in $x$-$y$ plane.}
\label{fig:Projection }
\end{figure}

\noindent
\begin{equation}\label{D12.1000}
\bullet \text{A definite integral }\int_{\theta=0}^{\frac{\pi}{2}}\sin^\alpha {\theta}\cos^\beta {\theta}~d\theta=\frac{\Gamma\left(\frac{\alpha+1}{2}\right) \Gamma\left(\frac{\beta+1}{2}\right) }{2\Gamma\left(\frac{\alpha+\beta+2}{2}\right) },~~~~~~~~~~~~~~~~~~~~~~~~~~~~~~~~~~~~~~~~~~~~~~~~
\end{equation}
where $\mathfrak{R}(\alpha)>-1,~\mathfrak{R}(\beta)>-1$.\\

Motivated by the work of Andrews\cite{And1985,And1998}, Bakshi {\em et al.}\cite{Bakshi1979}, Burchnall {\em et al.}\cite{Burch1941} and others \cite{Kampe26,ChhabRusia1980,Copson,Deshpande1981,Gradshteyn,Haisrimar92115,Humb20430,Humb2196,Humb217396,Kampe1921,MagOberSoni1966,MathaiHaubold,MathaiSaxenaHaubold,Qureshinaveddil2019498,Qureshinaveddil201952,Qureshinaveddil2020,Qureshinavediftara2020,Srivastava1967,Srichoi2011,SriGupta,SriPanda1976,Szego,Wright1935,Wright1940}, we evaluated some important definite integrals $\int_{\theta=-\pi}^{\pi}\left(\frac{\cos^2{\theta}}{\beta^2}+\frac{\sin^2{\theta}}{\lambda^2}\right)^sd\theta$ and $\int_{r=0}^{1} \frac{r^{2s+1}}{(1-r^2)^s}dr$ with suitable convergence conditions in {\bf section \ref{Z12.2}}, by using Mellin-Barnes type contour integral representation of binomial function$~_1F_0(z)$, Meijer's G-function, classical Beta function of two variables and series manipulation technique. These integrals are useful and help us in the derivation of closed form for the exact curved surface area of a hemiellipsoid. In {\bf section \ref{Z12.3}}, we derive the closed form of the exact curved surface area of a hemiellipsoid by using series manipulation technique, Mellin Barnes type contour integral representations of generalized hypergeometric function$~_pF_q(z)$, Meijer's $G$-function and analytic continuation formula for Gauss function in terms of Appell's function of first kind. In {\bf section \ref{Z12.5}}, we derive some special cases related to the total curved surface areas of ellipsoid, Prolate spheroid, Oblate spheroid and sphere. 

\section{\bf Evaluation of some useful definite integrals}\label{Z12.2}

\noindent
The following definite integrals hold true associated with suitable convergence conditions:
\begin{equation}\label{G12.10}
\text{\bf Theorem 1.}~\int_{\theta=-\pi}^{\pi}\left(\frac{\cos^2{\theta}}{\beta^2}+\frac{\sin^2{\theta}}{\lambda^2}\right)^s~d\theta =\frac{2\pi \lambda}{\beta^{1+2s}}~{}_2F_1\left[\begin{array}{ll}
\frac{1}{2},1+s;\\
&1-\frac{\lambda^2}{\beta^2}\\
~~~~~~~~1;\end{array}\right],~~
\end{equation}
where $\beta\geq \lambda>0$ and it is obvious that $0\leq(1-\frac{\lambda^2}{\beta^2})<1$.

\begin{equation}\label{G12.10000}
\text{\bf Theorem 2.}~\int_{\theta=-\pi}^{\pi}\left(\frac{\cos^2{\theta}}{\beta^2}+\frac{\sin^2{\theta}}{\lambda^2}\right)^s~d\theta =\frac{2\pi \beta}{\lambda^{1+2s}}~{}_2F_1\left[\begin{array}{ll}
\frac{1}{2},1+s;\\
&1-\frac{\beta^2}{\lambda^2}\\
~~~~~~~~1;\end{array}\right],~~
\end{equation}
where $\lambda\geq \beta>0$ and it is obvious that $0\leq(1-\frac{\beta^2}{\lambda^2})<1$.\\

\begin{equation}\label{l12.1000}
\text{\bf Theorem 3.}~\int_{r=0}^{1} \frac{r^{2s+1}}{(1-r^2)^s}~~dr=\frac{\Gamma(1+s)\Gamma(1-s)}{2},~~~~~~~~~~~~~~~~~~~~~~~~~~~~~~~~~~~~~~~~~~~~~~~
\end{equation}
where $|\mathfrak{R}(s)|<1.$ \\

\noindent
{\bf{Remark:}} The above formulas (\ref{G12.10}), (\ref{G12.10000}) and (\ref{l12.1000})are also verified numerically using Mathematica program.\\

\noindent
{\bf Independent demonstration of the assertions (\ref{G12.10}) and (\ref{G12.10000}) }
\begin{equation*}
\text{Suppose} ~I_1=\int_{\theta=-\pi}^{\pi}\left(\frac{\cos^2{\theta}}{\beta^2}+\frac{\sin^2{\theta}}{\lambda^2}\right)^s~d\theta.~~~~~~~~~~~~~~~~~~~~~~~~~~~~~
\end{equation*}
\begin{equation*}
=\frac{4}{\lambda^{2s}}\int_{\theta=0}^{\frac{\pi}{2}}\left( \sin^2 \theta\right)^s\left\lbrace 1+\frac{\lambda^2\cos^2{\theta}}{\beta^2\sin^2{\theta}}\right\rbrace ^s~d\theta 
\end{equation*}
\begin{equation}\label{G12.100}
=\frac{4}{\lambda^{2s}}\int_{\theta=0}^{\frac{\pi}{2}} \sin^{2s} \theta~{}_1F_0\left[\begin{array}{ll}
~~-s;\\
&\frac{-\lambda^2\cos^2{\theta}}{\beta^2\sin^2{\theta}}\\
~-~~;\end{array}\right]~d\theta. 
\end{equation}
Employing the contour integral (\ref{D12.100}) of ${}_{1}F_{0}(.)$, we get
\begin{equation}\label{G12.1005}
I_1=\frac{2}{\pi i \Gamma(-s)~\lambda^{2s}}\int_{\theta=0}^{\frac{\pi}{2}} \sin^{2s}\theta\left\lbrace \int_{\zeta=-i\infty}^{+i\infty}\Gamma(-\zeta)\Gamma(-s+\zeta)\left(\frac{\lambda^2\cos^2{\theta}}{\beta^2\sin^2{\theta}}\right)^\zeta~d\zeta \right\rbrace  ~d\theta.
\end{equation}
Interchanging the order of integration in double integral of (\ref{G12.1005}), we get
\begin{equation}
I_1=\frac{2}{\pi i \Gamma(-s)~\lambda^{2s}}\int_{\zeta=-i\infty}^{+i\infty}\Gamma(-\zeta)\Gamma(-s+\zeta)\left(\frac{\lambda^2}{\beta^2}\right)^\zeta\left\lbrace \int_{\theta=0}^{\frac{\pi}{2}}\sin^{2s-2\zeta}{\theta}~\cos^{2\zeta}{\theta}~d\theta\right\rbrace ~d\zeta.
\end{equation}
Using the integral formula (\ref{D12.1000}), we get
\begin{equation*}
I_1=\frac{1}{\pi i \Gamma(-s)\Gamma(1+s)~\lambda^{2s}}\int_{\zeta=-i\infty}^{+i\infty}\Gamma(-\zeta)\Gamma(-s+\zeta)\Gamma\left(\frac{1}{2}+s-\zeta\right)\Gamma\left(\frac{1}{2}+\zeta\right)\left(\frac{\lambda^2}{\beta^2}\right)^\zeta~d\zeta
\end{equation*}
\begin{equation}
=\frac{1}{\pi i \Gamma(-s)\Gamma(1+s)~\lambda^{2s}}\int_{\zeta=-i\infty}^{+i\infty}\Gamma(0-\zeta)\Gamma(1-(s+1)+\zeta)\Gamma\left(\frac{1}{2}+s-\zeta\right)\Gamma\left(1-\frac{1}{2}+\zeta\right)\left(\frac{\lambda^2}{\beta^2}\right)^\zeta~d\zeta.
\end{equation}
Applying the definition (\ref{E12.12}) of Meijer's $G$-function, we get
\begin{equation}\label{G12.101}
I_1=\frac{2}{\Gamma(-s)\Gamma(1+s)~\lambda^{2s} }~G^{2~2}_{2~2}\left(\frac{\lambda^2}{\beta^2}\left|   \begin{array}{ll}
s+1,\frac{1}{2};-\\
0,\frac{1}{2}+s;-  \end{array}\right. \right).~~~~~~~~~~~~~~~~~~~~~~~~~~~~~~~~~~~~
\end{equation}
Employing the conversion formula (\ref{F12.200}) in equation (\ref{G12.101}), and after further simplification, we arrive at the result (\ref{G12.10}).\\

\noindent
The proof of the result (\ref{G12.10000}) follows the same steps as in the proof of (\ref{G12.10}). So we omit the details here.\\

\noindent
{\bf Independent demonstration of the assertion (\ref{l12.1000}) }
\begin{equation*}
\text{Suppose} ~I_2=\int_{r=0}^{1} \frac{r^{2s+1}}{(1-r^2)^s}~~dr=\int_{0}^{1}{r^{2s+1}(1-r^2)^{-s}}~dr.~~~~~~~~~~~~~~~~~~~~~~~~~~~`
\end{equation*}
Put $r^2=t$, therefore $dr=\frac{1}{2\sqrt{t}}~dt$, we get
\begin{equation*}
~I_2=\frac{1}{2}\int_{0}^{1}t^{s}~{(1-t)^{-s}}~~dt;~~|\mathfrak{R}(s)|<1.~~~~~~~~~~~~~~~~~~~~~~~~~~~~~~~~~~~~~~~~~
\end{equation*}
Applying the definition of classical Beta function of two variables, we arrive at the result (\ref{l12.1000}).\\
\noindent
Throughout the discussion in the next sections, we are assuming that $a\geq b\geq c>0$.
\section{\bf Closed forms for curved surface area of a hemiellipsoid }\label{Z12.3}
\noindent
{\bf Theorem 4.} The curved surface area of hemiellipsoid (whose axis is positive direction of $z$-axis) i.e, $z=c\left(1- \frac{x^2}{a^2}-\frac{y^2}{b^2}\right)^{\frac{1}{2}};~a\geq b\geq c>0$, lying above $x$-$y$ plane (i.e, $z=0$) 
 is given by:\\
\noindent
\begin{equation}\label{p12.1}
\hat{S_1}=
=\left( \frac{2\pi b^2c^2 }{a^2}\right)  ~F_1\left[\begin{array}{ll}
2;~\frac{1}{2},~\frac{1}{2};~\frac{3}{2};~1-\frac{b^2}{a^2},1-\frac{c^2}{a^2}\end{array}\right],~~~~~~~~~~~~~~~~~~~~~~
\end{equation}
where $|1-\frac{b^2}{a^2}|<1~\text{or}~0<\frac{b^2}{a^2}<2$ and $|1-\frac{c^2}{a^2}|<1~\text{or}~0<\frac{c^2}{a^2}<2$.\\

\noindent
{\bf Theorem 5.} The curved surface area of hemiellipsoid (whose axis is positive direction of $y$-axis) i.e, $y=b\left(1- \frac{x^2}{a^2}-\frac{z^2}{c^2}\right)^{\frac{1}{2}};~a\geq b\geq c> 0$, lying above $x$-$z$ plane (i.e, $y=0$)
 is given by (\ref{p12.1}).\\
\noindent
{\bf Theorem 6.}~ The curved surface area of hemiellipsoid (whose axis is positive direction of $x$-axis) i.e, $x=a\left(1-\frac{y^2}{b^2}- \frac{z^2}{c^2}\right)^{\frac{1}{2}};~a\geq b\geq c>0$, lying above $y$-$z$ plane (i.e, $x=0$). Here two cases arise:\\

\noindent
{\bf Case I.} When $a>b $ such that $|1-\frac{a^2}{b^2}|<1$ i.e, $0<\frac{a^2}{b^2}<2$, then curved surface area of hemiellipsoid $x=a\left(1-\frac{y^2}{b^2}- \frac{z^2}{c^2}\right)^{\frac{1}{2}}$ is given by
\begin{equation}\label{t12.1}
\hat{S_2}=
=\left( \frac{2\pi a^2c^2 }{b^2}\right)  ~F_1\left[\begin{array}{ll}
2;~\frac{1}{2},~\frac{1}{2};~\frac{3}{2};~1-\frac{c^2}{b^2},1-\frac{a^2}{b^2}\end{array}\right],~~~~~~~~~~~~~~~~~~~~~~
\end{equation}
where $|1-\frac{a^2}{b^2}|<1~\text{or}~0<\frac{a^2}{b^2}<2$ and $|1-\frac{c^2}{b^2}|<1~\text{or}~0<\frac{c^2}{b^2}<2$.\\

\noindent
{\bf Case II.} When $a>b$ such that $|1-\frac{a^2}{b^2}|>1$ i.e, $2<\frac{a^2}{b^2}<\infty$, then curved surface area of hemiellipsoid $x=a\left(1-\frac{y^2}{b^2}- \frac{z^2}{c^2}\right)^{\frac{1}{2}}$ is given by

\begin{equation*}
\hat{S_3}=\left( \frac{\pi^2 a^2c^2 }{2b\sqrt{(a^2-b^2)}}\right)
~F_1\left[\begin{array}{ll}
\frac{1}{2};~\frac{3}{2},~\frac{1}{2};~1;~\frac{b^2-c^2}{b^2},\frac{b^2-c^2}{b^2-a^2}\end{array}\right]-~~~~~~~~~~~~~~~~~~~~~~~~~~
\end{equation*}
\begin{equation}\label{z12.3}
-\left( \frac{2\pi a^2b^2c^2 }{3(a^2-b^2)^2}\right)
~F_1\left[\begin{array}{ll}
2;~\frac{3}{2},~\frac{1}{2};~\frac{5}{2};~\frac{b^2}{b^2-a^2},\frac{b^2-c^2}{b^2-a^2}\end{array}\right],~~~~~~~~~~~
\end{equation}
where $|1-\frac{c^2}{b^2}|<1~\text{or}~0<\frac{c^2}{b^2}<2;~|1-\frac{a^2}{b^2}|>1~\text{or}~|\frac{b^2}{b^2-a^2}|<1~\text{or}~2<\frac{a^2}{b^2}<\infty$ and $|\frac{b^2-c^2}{b^2-a^2}|<1$.\\

\noindent
{\bf Remark:} The above formulas (\ref{p12.1}) to (\ref{z12.3}) are  equivalent and are verified numerically through {\em Mathematica program}.\\

\noindent
{\bf Demonstration of the assertion (\ref{p12.1})}\\
 \noindent
Equation of an ellipsoid is given by
\begin{equation}\label{l12.100}
\frac{x^2}{a^2}+\frac{y^2}{b^2}+\frac{z^2}{c^2}=1;~~a\geq b\geq c>0~~~~~~~~~~~~~~~
\end{equation}
\begin{figure}[ht!]
\noindent\begin{minipage}{0.40\textwidth}
\includegraphics[width=\linewidth]{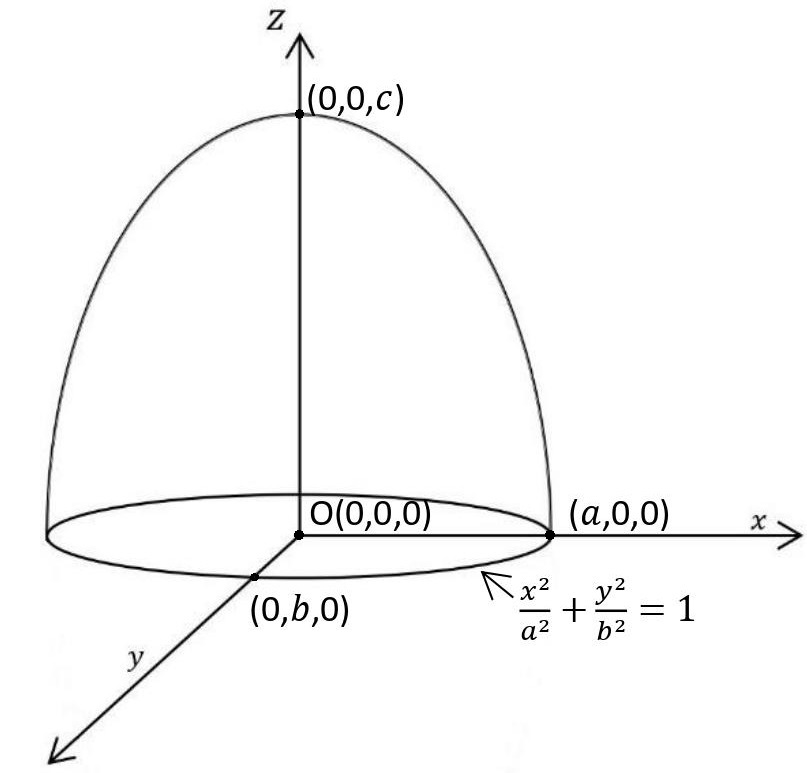}
\end{minipage}%
\caption[short]{Hemiellipsoid whose axis of symmetry is $z$- axis.}
\label{fig:Projection }
\end{figure}

\noindent
Therefore the equation of a hemiellipsoid will be
\begin{equation}\label{l12.101}
z=c\left(1- \frac{x^2}{a^2}-\frac{y^2}{b^2}\right)^{\frac{1}{2}}; ~c>0,~~~\text{taking positive sign of the square root for hemiellipsoid.}~~~~~~~~~~~~~~~~~~~~~~~~~~~~~~~~~~
\end{equation}
\begin{equation}
\text{Now~}~~~~\frac{\partial z}{\partial x}=\frac{-cx}{a^2\sqrt{\left(1- \frac{x^2}{a^2}-\frac{y^2}{b^2}\right)}}, ~~~~~~~~~~~~~~~~~~~~~~~~~~~~~~~~~~~~~~~~~~~
\end{equation}
\begin{equation}
\frac{\partial z}{\partial y}=\frac{-cy}{b^2\sqrt{\left(1- \frac{x^2}{a^2}-\frac{y^2}{b^2}\right)}}. ~~~~~~~~~~~~~~~~~~~~~~~~~~~~~~~~~~~~~~~~~~~
\end{equation}

\noindent
Substitute the values of $\frac{\partial z}{\partial x}$ and $\frac{\partial z}{\partial y}$ in equation (\ref{a12.9}).
Therefore curved surface area of hemiellipsoid will be
\begin{equation}
\hat{S_1}=\underbrace{\int\int}_{\stackrel{\text{over the area of the ellipse}}{{\frac{x^2}{a^2}+\frac{y^2}{b^2}=1}}}\sqrt{\left\lbrace 1+\frac{c^2x^2}{a^4\left(1- \frac{x^2}{a^2}-\frac{y^2}{b^2}\right)}+\frac{c^2y^2}{b^4\left(1-\frac{x^2}{a^2}-\frac{y^2}{b^2}\right)}\right\rbrace }~dx ~dy.
\end{equation}
Put $x=aX,~y=bY$, therefore, 
\begin{equation}
\hat{S_1}=ab~\underbrace{\int\int}_{\stackrel{\text{over the area of the circle}}{X^2+Y^2=1}}\sqrt{\left\lbrace 1+\frac{c^2X^2}{a^2\left(1- X^2-Y^2\right)}+\frac{c^2Y^2}{b^2\left(1- X^2-Y^2\right)}\right\rbrace  }~dX ~dY
\end{equation}
\begin{equation}
=ab~\underbrace{\int\int}_{\stackrel{\text{over the area of the circle}}{X^2+Y^2=1}}\sqrt{\left\lbrace 1+\frac{c^2}{\left(1- X^2-Y^2\right) }\left(  \frac{X^2}{a^2}+\frac{Y^2}{b^2}\right) \right\rbrace}~dX~dY.
\end{equation}

When $X=r\cos \theta,~Y=r\sin \theta$, then $dX~dY=rdr~d\theta$.
\begin{equation}\label{z12.102}
\text{Therefore}~ \hat{S_1}=ab~\int_{\theta=-\pi}^{\pi}\int_{r=0}^{1}\left\lbrace 1+\frac{c^2r^2}{(1-r^2)}\left( \frac{\cos^2\theta}{a^2}+\frac{\sin^2\theta}{b^2}\right)\right\rbrace^{\frac{1}{2}} ~rdr ~d\theta.
\end{equation}
{\bf Remark: }Since we have no standard formula of definite/indefinite integrals in the literature of integral calculus for the integration with respect to $"r"$ and $"\theta"$ in double integral (\ref{z12.102}). Therefore we can solve such integrals exactly through hypergeometric function approach.

\begin{equation}\label{l12.102}
\text{Therefore}~ \hat{S_1}=ab~\int_{\theta=-\pi}^{\pi} \int_{r=0}^{1}~_{1}F_{0}\left[\begin{array}{ll}
\frac{-1}{2};\\
& \frac{-c^2r^2}{(1-r^2)}\left( \frac{\cos^2\theta}{a^2}+\frac{\sin^2\theta}{b^2}\right) \\
~-;\end{array}\right]~rdr ~d\theta.
\end{equation}
Since there is uncertainty about the argument of$~_1F_0$ in equation (\ref{l12.102}), because the argument of $~_1F_0$ in equation (\ref{l12.102}) may be greater than 1. Therefore applying contour integral (\ref{D12.100}) of ${}_{1}F_{0}(.)$ in equation (\ref{l12.102}), we get

\begin{equation}\label{l12.30}
\hat{S_1}=ab~\int_{\theta=-\pi}^{\pi}\int_{r=0}^{1}\left[\frac{1}{(2\pi i)\Gamma(\frac{-1}{2})} ~\int_{s=-i\infty}^{+i\infty}~\Gamma(-s)\Gamma\left( \frac{-1}{2}+s\right) \left\lbrace \frac{c^2r^2}{(1-r^2)}\left( \frac{\cos^2\theta}{a^2}+\frac{\sin^2\theta}{b^2}\right) \right\rbrace^sds \right] ~rdr ~d\theta,
\end{equation}
where $|\arg\left\lbrace c^2\left( \frac{\cos^2\theta}{a^2}+\frac{\sin^2\theta}{b^2}\right) \right\rbrace|<\pi$ and $i=\sqrt{(-1)}.$\\

\noindent
Interchanging the order of integration in double integral of (\ref{l12.30}), we get
\begin{equation*}
\hat{S_1}=\frac{-ab}{4\pi\sqrt{\pi}~i}~\int_{s=-i\infty}^{+i\infty}~\Gamma(-s)\Gamma\left( \frac{-1}{2}+s\right)c^{2s}~\left\lbrace \int_{\theta=-\pi}^{{\pi}}\left( \frac{\cos^2\theta}{a^2}+\frac{\sin^2\theta}{b^2}\right)^s~d\theta \right\rbrace \times
\end{equation*}
\begin{equation}\label{l12.103}
\times~\left\lbrace \int_{r=0}^{1} \frac{r^{2s+1}}{(1-r^2)^s}~~dr\right\rbrace  ~ds.
\end{equation}
Since $a>b$ then employing the useful integrals (\ref{G12.10}) and (\ref{l12.1000}) in (\ref{l12.103}), we get
\begin{equation}
\hat{S_1}=\frac{-b^2}{4\sqrt{\pi}~i}~\int_{s=-i\infty}^{+i\infty}~ \Gamma(-s)\Gamma(1-s)\Gamma(1+s)\Gamma\left( \frac{-1}{2}+s\right)\frac{{c} ^{2s}}{a^{2s}}~{}_{2}F_{1}\left[\begin{array}{ll}
\frac{1}{2},1+s;\\
& 1-\frac{b^2}{a^2}\\
~~~~~~~~1;\end{array}\right]~ds, 
\end{equation}
where $|1-\frac{b^2}{a^2}|<1.$
\begin{equation}
\text{Therefore}~\hat{S_1}=\frac{-b^2}{4\sqrt{\pi}~i}~\int_{s=-i\infty}^{+i\infty} \Gamma(-s)\Gamma(1-s)\Gamma(1+s)\Gamma\left( \frac{-1}{2}+s\right)\frac{{c} ^{2s}}{a^{2s}}~\sum_{m=0}^{\infty}\frac{(\frac{1}{2})_m(1+s)_m(1-\frac{b^2}{a^2})^m}{(1)_m~m!}~ ds. 
\end{equation}

\begin{equation}
=\frac{-b^2}{4\sqrt{\pi}~i}\sum_{m=0}^{\infty}\frac{(\frac{1}{2})_m(1-\frac{b^2}{a^2})^m}{(1)_m~m!}~\int_{s=-i\infty}^{+i\infty} \Gamma(-s)\Gamma(1-s)\Gamma(1+s+m)\Gamma\left( \frac{-1}{2}+s\right)\left( \frac{{c^2} }{a^{2}}\right)^s ~ds. 
\end{equation}
Applying the definition (\ref{E12.12}) of Meijer's $G$-function, we get
\begin{equation}\label{l12.104}
\hat{S_1}=\frac{-b^2\sqrt{\pi}}{2}\sum_{m=0}^{\infty}\frac{(\frac{1}{2})_m(1-\frac{b^2}{a^2})^m}{(1)_m~m!}~G^{2~2}_{2~2}\left(\frac{c^2}{a^2}\left|   \begin{array}{ll}
-m,\frac{3}{2};-\\
~~~0,1;-  \end{array}\right. \right).~~~~~~~~~~~~~~~~~~~
\end{equation}
Employing the conversion formula (\ref{F12.200}) in equation (\ref{l12.104}), we get 
\begin{equation}
\hat{S_1}=\frac{2\pi b^2c^2}{a^2}\sum_{m=0}^{\infty}\frac{(\frac{1}{2})_m~(2)_m~(1-\frac{b^2}{a^2})^m}{(\frac{3}{2})_m~m!}~{}_{2}F_{1}\left[\begin{array}{ll}
\frac{1}{2},2+m;\\
& 1-\frac{c^2}{a^2}\\
~~\frac{3}{2}+m;\end{array}\right];~\left| 1-\frac{c^2}{a^2}\right| <1
\end{equation}
\begin{equation}\label{t12.105}
=\frac{2\pi b^2c^2}{a^2}\sum_{m=0}^{\infty}\frac{(\frac{1}{2})_m~(2)_m~(1-\frac{b^2}{a^2})^m}{(\frac{3}{2})_m~m!}~\sum_{k=0}^{\infty}\frac{(\frac{1}{2})_k~(2+m)_k~(1-\frac{c^2}{a^2})^k}{(\frac{3}{2}+m)_k~k!}.~~~~~~~~~
\end{equation}
Expressing the result (\ref{t12.105}) in terms of Appell's double hypergeometric function $F_1$, we get the result (\ref{p12.1}).\\

\noindent
{\bf Demonstration of the assertions (\ref{t12.1}) and (\ref{z12.3})}\\
 \noindent
The equation of a hemiellipsoid lying above $y$-$z$ plane is given by
\begin{equation}\label{t12.10}
x=a\left(1-\frac{y^2}{b^2}- \frac{z^2}{c^2}\right)^{\frac{1}{2}}; ~a\geq b\geq c>0,~\text{taking positive sign of the square root for hemiellipsoid.}~~~~~~~~~~~~~~~~~~~~~~~~~~~~~~~~~~
\end{equation}
\begin{equation}
\text{Now~}~~~~\frac{\partial x}{\partial y}=\frac{-ay}{b^2\sqrt{\left(1-\frac{y^2}{b^2}- \frac{z^2}{c^2}\right)}}, ~~~~~~~~~~~~~~~~~~~~~~~~~~~~~~~~~~~~~~~~~~~
\end{equation}
\begin{equation}
\frac{\partial x}{\partial z}=\frac{-az}{c^2\sqrt{\left(1-\frac{y^2}{b^2} -\frac{z^2}{c^2}\right)}}. ~~~~~~~~~~~~~~~~~~~~~~~~~~~~~~~~~~~~~~~~~~~
\end{equation}
\begin{figure}[ht!]
\noindent\begin{minipage}{0.45\textwidth}
\includegraphics[width=\linewidth]{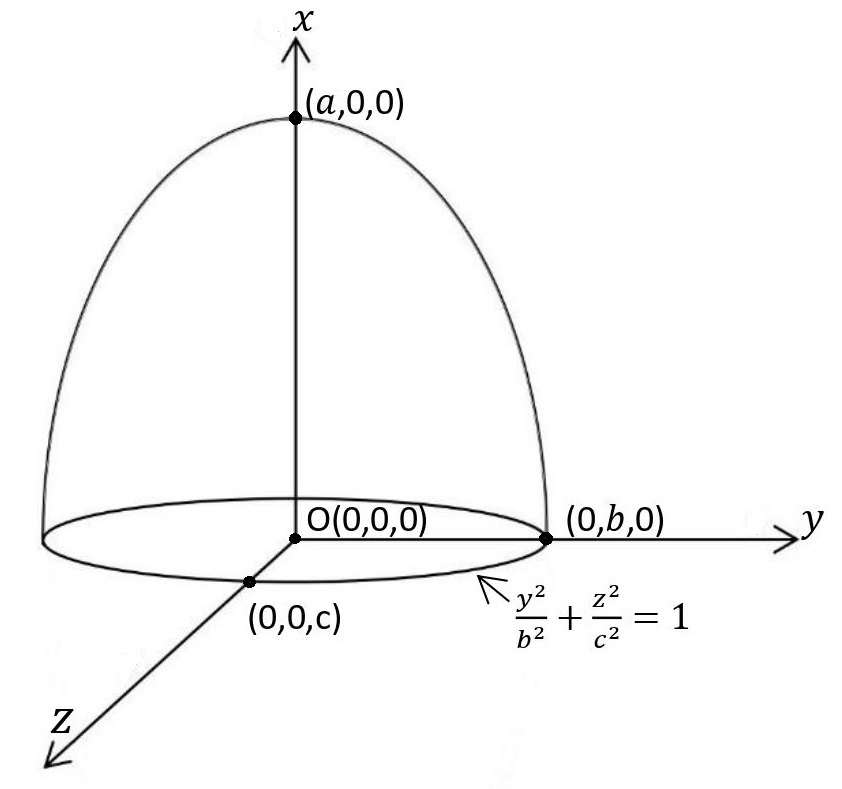}
\end{minipage}%
\caption[short]{Hemiellipsoid whose axis of symmetry is $y$- axis.}
\label{fig:Projection }
\end{figure}
\noindent
Substitute the values of $\frac{\partial x}{\partial y}$ and $\frac{\partial x}{\partial z}$ in equation (\ref{t12.9}).
Therefore curved surface area of hemiellipsoid will be
\begin{equation}
\hat{S_2}=\underbrace{\int\int}_{\stackrel{\text{over the area of the ellipse}}{{\frac{y^2}{b^2}+\frac{z^2}{c^2}=1}}}\sqrt{\left\lbrace 1+\frac{a^2y^2}{b^4\left(1- \frac{y^2}{b^2}-\frac{z^2}{c^2}\right)}+\frac{a^2z^2}{c^4\left(1-\frac{y^2}{b^2}-\frac{z^2}{c^2}\right)}\right\rbrace }~dy ~dz.
\end{equation}
Put $y=bY,~z=cZ$, therefore, 
\begin{equation}
\hat{S_2}=bc~\underbrace{\int\int}_{\stackrel{\text{over the area of the circle}}{Y^2+Z^2=1}}\sqrt{\left\lbrace 1+\frac{a^2Y^2}{b^2\left(1-Y^2-Z^2\right)}+\frac{a^2Z^2}{c^2\left(1- Y^2-Z^2\right)}\right\rbrace  }~dY ~dZ
\end{equation}
\begin{equation}
=bc~\underbrace{\int\int}_{\stackrel{\text{over the area of the circle}}{Y^2+Z^2=1}}\sqrt{\left\lbrace 1+\frac{a^2}{\left(1- Y^2-Z^2\right) }\left(  \frac{Y^2}{b^2}+\frac{Z^2}{c^2}\right) \right\rbrace}~dY~dZ.
\end{equation}

When $Y=r\cos \theta,~Z=r\sin \theta$, then $dY~dZ=rdr~d\theta$.
\begin{equation}\label{z12.11}
\text{Therefore}~ \hat{S_2}=bc~\int_{\theta=-\pi}^{\pi}\int_{r=0}^{1}\left\lbrace 1+\frac{a^2r^2}{(1-r^2)}\left( \frac{\cos^2\theta}{b^2}+\frac{\sin^2\theta}{c^2}\right)\right\rbrace^{\frac{1}{2}} ~rdr ~d\theta.
\end{equation}
{\bf Remark: }Since we have no standard formula of definite/indefinite integrals in the literature of integral calculus for the integration with respect to $"r"$ and $"\theta"$ in double integral (\ref{z12.11}). Therefore we can solve such integrals exactly through hypergeometric function approach.

\begin{equation}\label{t12.11}
\text{Therefore}~ \hat{S_2}=bc~\int_{\theta=-\pi}^{\pi} \int_{r=0}^{1}~_{1}F_{0}\left[\begin{array}{ll}
\frac{-1}{2};\\
& \frac{-a^2r^2}{(1-r^2)}\left( \frac{\cos^2\theta}{b^2}+\frac{\sin^2\theta}{c^2}\right) \\
~-;\end{array}\right]~rdr ~d\theta.
\end{equation}
Since there is uncertainty about the argument of$~_1F_0$ in equation (\ref{t12.11}), because the argument of $~_1F_0$ in equation (\ref{t12.11}) may be greater than 1. Therefore applying contour integral (\ref{D12.100}) of ${}_{1}F_{0}(.)$ in equation (\ref{t12.11}), we get
\begin{equation}\label{t12.12}
\hat{S_2}=bc~\int_{\theta=-\pi}^{\pi}\int_{r=0}^{1}\left[\frac{1}{(2\pi i)\Gamma(\frac{-1}{2})} ~\int_{s=-i\infty}^{+i\infty}~\Gamma(-s)\Gamma\left( \frac{-1}{2}+s\right) \left\lbrace \frac{a^2r^2}{(1-r^2)}\left( \frac{\cos^2\theta}{b^2}+\frac{\sin^2\theta}{c^2}\right) \right\rbrace^sds \right] ~rdr ~d\theta,
\end{equation}
where $|\arg\left\lbrace a^2\left( \frac{\cos^2\theta}{b^2}+\frac{\sin^2\theta}{c^2}\right) \right\rbrace|<\pi$ and $i=\sqrt{(-1)}.$\\

\noindent
Interchanging the order of integration in double integral of (\ref{t12.12}), we get
\begin{equation*}
\hat{S_2}=\frac{-bc}{4\pi\sqrt{\pi}~i}~\int_{s=-i\infty}^{+i\infty}~\Gamma(-s)\Gamma\left( \frac{-1}{2}+s\right)a^{2s}~\left\lbrace \int_{\theta=-\pi}^{{\pi}}\left( \frac{\cos^2\theta}{b^2}+\frac{\sin^2\theta}{c^2}\right)^s~d\theta \right\rbrace \times
\end{equation*}
\begin{equation}\label{t12.13}
\times~\left\lbrace \int_{r=0}^{1} \frac{r^{2s+1}}{(1-r^2)^s}~~dr\right\rbrace  ~ds.
\end{equation}
Since $b>c$ then employing the useful integrals (\ref{G12.10}) and (\ref{l12.1000}) in (\ref{t12.13}), we get
\begin{equation}
\hat{S_2}=\frac{-c^2}{4\sqrt{\pi}~i}~\int_{s=-i\infty}^{+i\infty}~ \Gamma(-s)\Gamma(1-s)\Gamma(1+s)\Gamma\left( \frac{-1}{2}+s\right)\frac{{a} ^{2s}}{b^{2s}}~{}_{2}F_{1}\left[\begin{array}{ll}
\frac{1}{2},1+s;\\
& 1-\frac{c^2}{b^2}\\
~~~~~~~~1;\end{array}\right]~ds, 
\end{equation}
where $|1-\frac{c^2}{b^2}|<1.$
\begin{equation}
\text{Therefore}~\hat{S_2}=\frac{-c^2}{4\sqrt{\pi}~i}~\int_{s=-i\infty}^{+i\infty} \Gamma(-s)\Gamma(1-s)\Gamma(1+s)\Gamma\left( \frac{-1}{2}+s\right)\frac{{a} ^{2s}}{b^{2s}}~\sum_{m=0}^{\infty}\frac{(\frac{1}{2})_m(1+s)_m(1-\frac{c^2}{b^2})^m}{(1)_m~m!}~ ds. 
\end{equation}

\begin{equation}
=\frac{-c^2}{4\sqrt{\pi}~i}\sum_{m=0}^{\infty}\frac{(\frac{1}{2})_m(1-\frac{c^2}{b^2})^m}{(1)_m~m!}~\int_{s=-i\infty}^{+i\infty} \Gamma(-s)\Gamma(1-s)\Gamma(1+s+m)\Gamma\left( \frac{-1}{2}+s\right)\left( \frac{{a^2} }{b^{2}}\right)^s ~ds. 
\end{equation}
Applying the definition (\ref{E12.12}) of Meijer's $G$-function, we get
\begin{equation}\label{t12.14}
\hat{S_2}=\frac{-c^2\sqrt{\pi}}{2}\sum_{m=0}^{\infty}\frac{(\frac{1}{2})_m(1-\frac{c^2}{b^2})^m}{(1)_m~m!}~G^{2~2}_{2~2}\left(\frac{a^2}{b^2}\left|   \begin{array}{ll}
-m,\frac{3}{2};-\\
~~~0,1;-  \end{array}\right. \right).~~~~~~~~~~~~~~~~~~~
\end{equation}
Employing the conversion formula (\ref{F12.200}) in equation (\ref{t12.14}), we get 
\begin{equation}
\hat{S_2}=\frac{2\pi a^2c^2}{b^2}\sum_{m=0}^{\infty}\frac{(\frac{1}{2})_m~(2)_m~(1-\frac{c^2}{b^2})^m}{(\frac{3}{2})_m~m!}~{}_{2}F_{1}\left[\begin{array}{ll}
\frac{1}{2},2+m;\\
& 1-\frac{a^2}{b^2}\\
~~\frac{3}{2}+m;\end{array}\right];~|1-\frac{a^2}{b^2}|<1
\end{equation}
\begin{equation}\label{t12.15}
=\frac{2\pi a^2c^2}{b^2}\sum_{m=0}^{\infty}\frac{(\frac{1}{2})_m~(2)_m~(1-\frac{c^2}{b^2})^m}{(\frac{3}{2})_m~m!}~\sum_{k=0}^{\infty}\frac{(\frac{1}{2})_k~(2+m)_k~(1-\frac{a^2}{b^2})^k}{(\frac{3}{2}+m)_k~k!}.~~~~~~~~~
\end{equation}

\noindent
{\bf Case I} (Corollary 1)\\
When $~a\geq b\geq c>0$ and $|1-\frac{a^2}{b^2}|<1~\text{or}~0<\frac{a^2}{b^2}<2$, then $|1-\frac{c^2}{b^2}|<1~\text{or}~0<\frac{c^2}{b^2}<2$ is always possible, therefore expressing the result (\ref{t12.15}) in terms of Appell's double hypergeometric function $F_1$, we get the result (\ref{t12.1}) .\\

\noindent
{\bf Case II} (Corollary 2)\\
When $~a\geq b>c>0$ and $|1-\frac{a^2}{b^2}|>1~\text{or}~2<\frac{a^2}{b^2}<\infty$, therefore from the result (\ref{t12.15}), we have
\begin{equation}\label{t12.16}
~\hat{S_3}=\left( \frac{2\pi a^2c^2 }{b^2}\right) \sum_{m=0}^{\infty}\frac{\left( 2\right)_{m} \left( \frac{1}{2}\right)_m \left( 1-\frac{c^2}{b^2}\right)^m }{\left( \frac{3}{2}\right)_m ~m!} ~{}_2F_1\left[\begin{array}{ll}
\frac{1}{2},~2+m;\\
&1-\frac{a^2}{b^2}\\
~~~~\frac{3}{2}+m;\end{array}\right].
\end{equation}
When $|1-\frac{a^2}{b^2}|>1~\text{or}~2<\frac{a^2}{b^2}<\infty$, then applying analytic continuation formula (\ref{G7.12}) for Gauss' function $_{2}F_1$ in equation (\ref{t12.16}), we get

\begin{equation*}
\hat{S_3}=\left( \frac{2\pi a^2c^2}{b^2}\right) \sum_{m=0}^{\infty}\frac{(2)_m\left( \frac{1}{2}\right)_m \left( 1-\frac{c^2}{b^2}\right)^m }{\left( \frac{3}{2}\right)_m ~m!}\left\lbrace \frac{\Gamma(\frac{3}{2}+m)\Gamma(\frac{-3}{2}-m)}{\Gamma(\frac{1}{2})\Gamma(\frac{-1}{2})}\left( \frac{a^2}{b^2}-1\right)^{-(2+m)}\times  \right.
\end{equation*}
\begin{equation}
\left.\times~{}_2F_1\left[\begin{array}{ll}
\frac{3}{2},2+m;\\
& \frac{b^2}{b^2-a^2}\\
~~\frac{5}{2}+m;\end{array}\right]
+\frac{\Gamma(\frac{3}{2}+m)\Gamma(\frac{3}{2}+m)}{\Gamma(1+m)\Gamma(2+m)}\left( \frac{a^2}{b^2}-1\right)^{-\frac{1}{2}}~{}_2F_1\left[\begin{array}{ll}
~~-m,~\frac{1}{2};\\
& \frac{b^2}{b^2-a^2}\\
-\frac{1}{2}-m;\end{array}\right]\right\rbrace 
\end{equation}
\begin{equation*}
=\left( \frac{2\pi a^2c^2}{3b^2}\right) \sum_{m=0}^{\infty}\frac{(2)_m\left( \frac{1}{2}\right)_m \left( 1-\frac{c^2}{b^2}\right)^m (-1)^{m+1}\left(\frac{b^2}{a^2-b^2} \right)^{m+2} }{\left( \frac{5}{2}\right)_m ~m!}~{}_2F_1\left[\begin{array}{ll}
\frac{3}{2},2+m;\\
& \frac{b^2}{b^2-a^2}\\
~~\frac{5}{2}+m;\end{array}\right]+
\end{equation*}
\begin{equation}
+\left( \frac{\pi^2 a^2c^2}{2b^2}\right) \sum_{m=0}^{\infty}\frac{\left( \frac{1}{2}\right)_m \left( \frac{3}{2}\right)_m\left( 1-\frac{c^2}{b^2}\right)^m \left(\frac{b^2}{a^2-b^2} \right)^{\frac{1}{2}} }{\left(1\right)_m ~m!}~{}_2F_1\left[\begin{array}{ll}
~~-m,~\frac{1}{2};\\
& \frac{b^2}{b^2-a^2}\\
-\frac{1}{2}-m;\end{array}\right]
\end{equation}
\begin{equation*}
=\left( \frac{\pi^2 a^2c^2 }{2b\sqrt{(a^2-b^2)}}\right) \sum_{m=0}^{\infty}\frac{\left( \frac{1}{2}\right)_m \left( \frac{3}{2}\right)_m \left( 1-\frac{c^2}{b^2}\right)^m }{\left( 1\right)_m~m!} ~{}_2F_1\left[\begin{array}{ll}
~~~-m,~\frac{1}{2};\\
&\frac{b^2}{b^2-a^2}\\
~-\frac{1}{2}-m~;\end{array}\right]-~~~~~~~~~~~~~~~~~~
\end{equation*}
\begin{equation}
-\left( \frac{2\pi a^2b^2c^2 }{3(a^2-b^2)^2}\right) \sum_{m=0}^{\infty}\frac{\left( 2\right)_m \left( \frac{1}{2}\right)_m \left( \frac{b^2-c^2}{b^2-a^2}\right)^m }{\left( \frac{5}{2}\right)_m ~m!} ~{}_2F_1\left[\begin{array}{ll}
2+m,~\frac{3}{2};\\
&\frac{b^2}{b^2-a^2}\\
~~~\frac{5}{2}+m~;\end{array}\right],~~~~
\end{equation}
where $|1-\frac{c^2}{b^2}|<1~\text{or}~0<\frac{c^2}{b^2}<2;~|1-\frac{a^2}{b^2}|>1~\text{or}~~|\frac{b^2}{b^2-a^2}|<1~\text{or}~2<\frac{a^2}{b^2}<\infty$ and $|\frac{b^2-c^2}{b^2-a^2}|<1$.

\begin{equation*}
\text{Therefore}~\hat{S_3}=\left( \frac{\pi^2 a^2c^2 }{2b\sqrt{(a^2-b^2)}}\right) \sum_{m=0}^{\infty}\frac{\left( \frac{1}{2}\right)_m \left( \frac{3}{2}\right)_m \left( 1-\frac{c^2}{b^2}\right)^m }{\left( 1\right)_m~m!}\sum_{r=0}^{m} \frac{(-m)_r\left(\frac{1}{2}\right)_r\left(\frac{b^2}{b^2-a^2} \right)^r  }{\left(-\frac{1}{2}-m\right)_r~r! }-~~~~~~~~~~~~~~~~~~
\end{equation*}
\begin{equation}
-\left( \frac{2\pi a^2b^2c^2 }{3(a^2-b^2)^2}\right) \sum_{m=0}^{\infty}\frac{\left( 2\right)_m \left( \frac{1}{2}\right)_m \left( \frac{b^2-c^2}{b^2-a^2}\right)^m }{\left( \frac{5}{2}\right)_m ~m!}\sum_{r=0}^{\infty}\frac{ 
(2+m)_r\left( \frac{3}{2}\right)_r\left( \frac{b^2}{b^2-a^2}\right)^r}{\\
\left( \frac{5}{2}+m\right)_r~r!}~~~~
\end{equation}
\begin{equation*}
=\left( \frac{\pi^2 a^2c^2 }{2b\sqrt{(a^2-b^2)}}\right) \sum_{m=0}^{\infty}\frac{\left( \frac{1}{2}\right)_m \left( \frac{3}{2}\right)_m \left( 1-\frac{c^2}{b^2}\right)^m }{~m!}\sum_{r=0}^{m} \frac{\left(\frac{1}{2}\right)_r\left(\frac{b^2}{a^2-b^2} \right)^r  }{\left(-\frac{1}{2}-m\right)_r~(m-r)!~r! }-~~~~~~~~~~~~~~~~~~
\end{equation*}
\begin{equation}\label{t12.17}
-\left( \frac{2\pi a^2b^2c^2 }{3(a^2-b^2)^2}\right) \sum_{m=0}^{\infty}\sum_{r=0}^{\infty}\frac{\left( 2\right)_{m+r} \left( \frac{1}{2}\right)_m \left( \frac{3}{2}\right)_r\left( \frac{b^2-c^2}{b^2-a^2}\right)^m \left( \frac{b^2}{b^2-a^2}\right)^r}{\left( \frac{5}{2}\right)_{m+r} ~m!~r!}~~~~
\end{equation}
Replacing $m$ by $m+r$ in the first series on right hand side of equation (\ref{t12.17}), we get

\begin{equation*}
\hat{S_3}=\left( \frac{\pi^2 a^2c^2 }{2b\sqrt{(a^2-b^2)}}\right) \sum_{m=0}^{\infty}\sum_{r=0}^{\infty}\frac{\left( \frac{1}{2}\right)_{m+r} \left( \frac{3}{2}\right)_{m+r}\left(\frac{1}{2}\right)_r \left( 1-\frac{c^2}{b^2}\right)^{m+r}\left(\frac{b^2}{a^2-b^2} \right)^r  }{(1)_{m+r}\left(-\frac{1}{2}-m-r\right)_r~m!~r!} \frac{ }{ }-~~~~~~~~~~~~~~~~~~
\end{equation*}
\begin{equation}
-\left( \frac{2\pi a^2b^2c^2 }{3(a^2-b^2)^2}\right) \sum_{m=0}^{\infty}\sum_{r=0}^{\infty}\frac{\left( 2\right)_{m+r} \left( \frac{1}{2}\right)_m \left( \frac{3}{2}\right)_r\left( \frac{b^2-c^2}{b^2-a^2}\right)^m \left( \frac{b^2}{b^2-a^2}\right)^r}{\left( \frac{5}{2}\right)_{m+r} ~m!~r!}~~~~
\end{equation}
\begin{equation*}
\text{Therefore}~\hat{S_3}=\left( \frac{\pi^2 a^2c^2 }{2b\sqrt{(a^2-b^2)}}\right) \sum_{m=0}^{\infty}\sum_{r=0}^{\infty}\frac{\left( \frac{1}{2}\right)_{m+r} \left( \frac{3}{2}\right)_{m}\left(\frac{1}{2}\right)_r \left( 1-\frac{c^2}{b^2}\right)^{m}\left(\frac{b^2-c^2}{b^2-a^2} \right)^r  }{(1)_{m+r}~m!~r!} \frac{ }{ }-~~~~~~~~~~~~~~~~~~
\end{equation*}
\begin{equation}\label{t12.18}
-\left( \frac{2\pi a^2b^2c^2 }{3(a^2-b^2)^2}\right) \sum_{m=0}^{\infty}\sum_{r=0}^{\infty}\frac{\left( 2\right)_{m+r} \left( \frac{1}{2}\right)_m \left( \frac{3}{2}\right)_r\left( \frac{b^2-c^2}{b^2-a^2}\right)^m \left( \frac{b^2}{b^2-a^2}\right)^r}{\left( \frac{5}{2}\right)_{m+r} ~m!~r!},~~~~
\end{equation}
where $|1-\frac{c^2}{b^2}|<1~\text{or}~0<\frac{c^2}{b^2}<2;~|1-\frac{a^2}{b^2}|>1~\text{or}~~|\frac{b^2}{b^2-a^2}|<1~\text{or}~2<\frac{a^2}{b^2}<\infty$ and $|\frac{b^2-c^2}{b^2-a^2}|<1$.\\

\noindent
On using the definition (\ref{C1.4}) of Appell's function of first kind $F_1$ in right hand side of (\ref{t12.18}), we arrive at the result (\ref{z12.3}).\\

\noindent
{\bf Remark:} Similarly the derivation of the formula for the curved surface area of a hemiellipsoid $y=b\left(1- \frac{x^2}{a^2}-\frac{z^2}{c^2}\right)^{\frac{1}{2}};~a\geq b\geq c>0$, whose axis is $y$-axis, lying above $x$-$z$ plane follows the same steps as above. So we omit the details here and we find the formula (\ref{p12.1}).

\section{\bf Some formulas for total curved surface areas of ellipsoid, Prolate spheroid and Oblate spheroid}\label{Z12.5}
\noindent
In this section, we discuss some special cases of the closed forms (\ref{p12.1}) to (\ref{z12.3}):\\

\noindent 
{\bf (i)} For the total curved surface area of complete ellipsoid $\frac{x^2}{a^2}+\frac{y^2}{b^2}+\frac{z^2}{c^2}=1;a\geq b\geq c>0$ whose axis is $z$- axis, multiply the closed form (\ref{p12.1}) by 2, we get 
\begin{equation}\label{l12.150}
\hat{S_1}=\left( \frac{4\pi b^2c^2 }{a^2}\right)  ~F_1\left[\begin{array}{ll}
2;~\frac{1}{2},~\frac{1}{2};~\frac{3}{2};~1-\frac{b^2}{a^2},1-\frac{c^2}{a^2}\end{array}\right];~\left|1-\frac{b^2}{a^2}\right|<1,\left|1-\frac{c^2}{a^2}\right|<1.~~~~~~~~~~~~~~~
\end{equation}

\noindent
{\bf (ii)} On substituting $b=a$ in equation (\ref{l12.150}), 
we get the total curved surface area of Oblate spheroid $\frac{x^2+y^2}{a^2}+\frac{z^2}{c^2}=1;a>c,$ formed by revolving the ellipse ($\frac{y^2}{b^2}+\frac{z^2}{c^2}=1$ lying in $y$-$z$ plane) or ($\frac{x^2}{a^2}+\frac{z^2}{c^2}=1$ lying in $x$-$z$ plane) about $z$-axis (i.e, minor axis)
\begin{equation}\label{l12.151}
\hat{S_1}=4\pi c^2 ~_2{F}_1\left[\begin{array}{lll}
2,\frac{1}{2};\\
& 1-\frac{c^2}{a^2}\\
~~~\frac{3}{2};\\\end{array}\right];~~0<\left(1-\frac{c^2}{a^2}\right)<2.~~~~~~~~~~
\end{equation}
Applying the formula (\ref{Z1.9}) in equation (\ref{l12.151}), we get the curved surface area of Oblate spheroid
\begin{equation}\label{X12.151}
\hat{S_1}=2\pi a^2 \left(1+\frac{c^2\tanh^{-1}\left( \sqrt{1-\frac{c^2}{a^2}}\right) }{a\sqrt{(a^2-c^2)}}\right),~~~~~~~~~~~~~~~~~~~~~~~~~
\end{equation}
\begin{equation}
\text{or}~\hat{S_1}=2\pi a^2 +\frac{\pi c^2}{\sqrt{\left( 1-\frac{c^2}{a^2}\right) }}\ln \left\lbrace \frac{1+ \sqrt{\left(1-\frac{c^2}{a^2}\right)} }{1- \sqrt{\left(1-\frac{c^2}{a^2}\right)}}\right\rbrace .~~~~~~~~~~~~~~~~~
\end{equation}

\noindent 
{\bf (iii)} For the total curved surface area of complete ellipsoid $\frac{x^2}{a^2}+\frac{y^2}{b^2}+\frac{z^2}{c^2}=1;a\geq b\geq c>0$ whose axis is $x$- axis, multiply the closed forms (\ref{t12.1}) and (\ref{z12.3}) by 2, we get 
\begin{equation}\label{L12.150}
\hat{S_2}=\left( \frac{4\pi a^2c^2 }{b^2}\right)  ~F_1\left[\begin{array}{ll}
2;~\frac{1}{2},~\frac{1}{2};~\frac{3}{2};~1-\frac{c^2}{b^2},1-\frac{a^2}{b^2}\end{array}\right],~~~~~~~~~~~~~~~~~~~~
\end{equation}
where $\left| 1-\frac{a^2}{b^2}\right| <1~\text{or}~0<\frac{a^2}{b^2}<2$ and $\left| 1-\frac{c^2}{b^2}\right| <1~\text{or}~0<\frac{c^2}{b^2}<2$.
\begin{equation*}
\text{and}~\hat{S_3}=\left( \frac{\pi^2 a^2c^2 }{b\sqrt{(a^2-b^2)}}\right)
~F_1\left[\begin{array}{ll}
\frac{1}{2};~\frac{3}{2},~\frac{1}{2};~1;~\frac{b^2-c^2}{b^2},\frac{b^2-c^2}{b^2-a^2}\end{array}\right]-~~~~~~~~~~~~~~~~~~~~~~~~~~
\end{equation*}
\begin{equation}\label{L12.3}
-\left( \frac{4\pi a^2b^2c^2 }{3(a^2-b^2)^2}\right)
~F_1\left[\begin{array}{ll}
2;~\frac{3}{2},~\frac{1}{2};~\frac{5}{2};~\frac{b^2}{b^2-a^2},\frac{b^2-c^2}{b^2-a^2}\end{array}\right],~~~~~~~~~~~
\end{equation}
where $|1-\frac{c^2}{b^2}|<1~\text{or}~0<\frac{c^2}{b^2}<2;~|1-\frac{a^2}{b^2}|>1~\text{or}~|\frac{b^2}{b^2-a^2}|<1~\text{or}~2<\frac{a^2}{b^2}<\infty$ and $|\frac{b^2-c^2}{b^2-a^2}|<1$.\\

\noindent
{\bf (iv)} Substitute $c=b$ in equation (\ref{L12.150}), we get the total curved surface area of Prolate spheroid $\frac{x^2}{a^2}+\frac{y^2+z^2}{b^2}=1;a>b$ such that $\left|1-\frac{a^2}{b^2}\right|<1$, formed by revolving the ellipse ($\frac{x^2}{a^2}+\frac{y^2}{b^2}=1$ lying in $x$-$y$ plane) or ($\frac{x^2}{a^2}+\frac{z^2}{c^2}=1$ lying in $x$-$z$ plane) about $x$-axis (i.e, major axis)
\begin{equation}\label{m12.151}
\hat{S_2}=4\pi a^2 ~_2{F}_1\left[\begin{array}{lll}
2,~\frac{1}{2};\\
& 1-\frac{a^2}{b^2}\\
~~~~\frac{3}{2};\\\end{array}\right],~~~~~~~~~~~~~~~~~~~~~~~~~~~~~~~~~~~~~~~~~~~~
\end{equation}
where $0<\frac{a^2}{b^2}<2.$\\
Applying the formula (\ref{Z1.9}) in equation (\ref{m12.151}), we get  the curved surface area of Prolate spheroid
\begin{equation}\label{m12.158}
\hat{S_2}=2\pi b^2 \left(1+\frac{a^2\tanh^{-1}\left( \sqrt{1-\frac{a^2}{b^2}}\right) }{b^2\sqrt{(1-\frac{a^2}{b^2})}}\right),~~~~~~~~~~~~~~~~~~~~~~~~~~~~~~~~~~~~~~~~~~~~~~~~~~~
\end{equation}
where $-\infty<\left(1-\frac{a^2}{b^2}\right)<1$.\\
or
\begin{equation}\label{P12.158}
\hat{S_2}=2\pi b^2 \left(1+\frac{a^2\tan^{-1}\left( \sqrt{\frac{a^2-b^2}{b^2}}\right) }{b\sqrt{(a^2-b^2)}}\right),~~~~~~~~~~~~~~~~~~~~~~~~~~~~~~~~~~~~~~~~~
\end{equation}
where $0\leq\left(\frac{a^2-b^2}{b^2}\right)<\infty.$

\begin{equation}
\text{or}~\hat{S_2}=2\pi b^2+\frac{2\pi ab}{\sqrt{\left( 1-\frac{b^2}{a^2}\right) }} \sin^{-1}\left( \sqrt{1-\frac{b^2}{a^2}}\right) .~~~~~~~~~~~~~~~~~~~~~~~~~
\end{equation}
\noindent
{\bf (v)} Substitute $c=b$ in equation (\ref{L12.3}), we get the total curved surface area of Prolate spheroid $\frac{y^2+z^2}{b^2}+\frac{x^2}{a^2}=1;a>b$ such that $\left|1-\frac{a^2}{b^2}\right|>1~\text{or}~2<\frac{a^2}{b^2}<\infty$, formed by revolving the ellipse ($\frac{x^2}{a^2}+\frac{y^2}{b^2}=1$ lying in $x$-$y$ plane) or ($\frac{x^2}{a^2}+\frac{z^2}{c^2}=1$ lying in $x$-$z$ plane) about $x$-axis (i.e, major axis)
\begin{equation}\label{p12.151}
\hat{S_3}=\frac{\pi^2 a^2 b}{\sqrt{(a^2-b^2)}}-\frac{4\pi a^2b^4}{3(a^2-b^2)^2} ~_2{F}_1\left[\begin{array}{lll}
2,~\frac{3}{2};\\
& \frac{b^2}{b^2-a^2}\\
~~~~\frac{5}{2};\\\end{array}\right],~~~~~~~~~~~~~~~~~~~~~~~~~~~~~~~
\end{equation}
where $0<\left(\frac{b^2}{b^2-a^2}\right)<2.$\\
Applying the formula (\ref{Z1.8}) in equation (\ref{p12.151}), we get 

\begin{equation}\label{Z12.151}
\hat{S_3}=\frac{\pi^2 a^2 b}{\sqrt{(a^2-b^2)}}-\frac{2\pi a^2 b^2}{(a^2-b^2)} \left(\frac{(b^2-a^2)}{a^2}+\frac{\tanh^{-1}\left( \sqrt{\frac{b^2}{b^2-a^2}}\right) }{\sqrt{(\frac{b^2}{b^2-a^2})}}\right),~~~~~~~~~~~~~~
\end{equation}
where $-\infty<\left(\frac{b^2}{b^2-a^2}\right)<1,\left(\frac{b^2}{b^2-a^2}\right)\neq 0.$\\
or
\begin{equation}\label{X12.155}
\hat{S_3}=2\pi b^2 \left(1+\frac{\pi a^2}{2b\sqrt{(a^2-b^2)}}-\frac{a^2\tan^{-1}\left( \sqrt{\frac{b^2}{a^2-b^2}}\right) }{b\sqrt{(a^2-b^2)}}\right),~~~~~~~~~~~~~~~~~~~~~~~~~~~~
\end{equation}
where $0\leq\left(\frac{b^2}{a^2-b^2}\right)<\infty$.\\

\noindent
{\bf Remark:} The above formulas (\ref{m12.158}) and (\ref{Z12.151}) are also valid in the interval $-\infty<\left(1-\frac{a^2}{b^2}\right)<1$ and $-\infty<\left(\frac{b^2}{b^2-a^2}\right)<1$, respectively. We have checked it using {\it Mathematica Program}.\\

\noindent
{\bf (vi)} Put $c=a$ and $b=a$ in equation (\ref{l12.150}), we get the total curved surface area of a sphere $x^2+y^2+z^2=a^2$ which is given by 
$\hat{S_1}={4\pi a^2}$.

\section{\bf Conclusion}
In this paper, we obtained the closed form for the exact curved surface area of a hemiellipsoid $z=c\left( 1-\frac{x^2}{a^2}-\frac{y^2}{b^2}\right)^{\frac{1}{2}}$ through hypergeometric function approach i.e, by using series rearrangement technique, Mellin-Barnes type contour integral representations of generalized hypergeometric function$~_pF_q(z)$, Meijer's $G$-function and analytic continuation formula for Gauss function; in terms of Appell's function of first kind. These formulas are neither available in the literature of mathematics nor found in any mathematical tables. Moreover, we also derived some special cases related to ellipsoid, Prolate spheroid, Oblate spheroid and sphere. We conclude that many formulas for the curved surface areas of other three dimensional figures can be derived in an analogous manner, using Mellin-Barnes contour integration. Moreover, the results deduced above (presumably new), have potential applications in the fields of applied mathematics, statistics and engineering sciences.\\

\noindent
{\bf Conflicts of interests:} The authors declare that there are no conflicts of interest.

\end{document}